\documentclass[12pt]{article}

\usepackage[T1]{fontenc}
\usepackage[utf8]{inputenc}
\usepackage{lmodern}
\usepackage{microtype}

\usepackage[left=1.05in,right=1.05in,top=1in,bottom=1in]{geometry}
\usepackage{amsmath,amssymb,amsfonts,amsthm,mathtools}
\usepackage{graphicx}
\usepackage{float}
\usepackage{tikz}
\usetikzlibrary{arrows.meta,calc,positioning,decorations.markings}
\usepackage{hyperref}
\usepackage[noabbrev]{cleveref}

\hypersetup{
  colorlinks=true,
  linkcolor=blue!50!black,
  citecolor=blue!60!black,
  urlcolor=blue!60!black
}

\numberwithin{equation}{section}

\theoremstyle{plain}
\newtheorem{theorem}{Theorem}[section]
\newtheorem{proposition}[theorem]{Proposition}

\theoremstyle{definition}

\theoremstyle{remark}
\newtheorem{remark}[theorem]{Remark}
\newtheorem{example}[theorem]{Example}

\title{Henneberg-Type Minimal Immersions in $\mathbb{R}^4$}
\author{Erhan G\"uler and Magdalena Toda}
\date{}

\begin{document}
\maketitle

\begin{abstract}
{\color{black}
We construct a family of higher-order Henneberg-type minimal surfaces in
$\mathbb{R}^4$ from explicit holomorphic data in a generalized
Weierstrass--Enneper representation. The construction retains the
polynomial reciprocal structure of the classical Henneberg data and yields
explicit holomorphic primitives and real parametrizations for the family
$\mathcal H_{(m,n)}$. We distinguish regular points from branch points,
analyze the role of the complex parameter $\lambda$, and show that the
real-$\lambda$ subfamily $\mathcal H_{(1,1)}$ lies in a three-dimensional
linear subspace. For $\mathcal H_{(1,1)}$, we compute the induced conformal
metric and Gaussian curvature, proving that the curvature is negative at
every regular point. An explicit parametrization is given for $\mathcal H_{(1,1)}$, together with three-dimensional coordinate
projections used solely for visualization. We also derive an integral-free
formulation in terms of a single holomorphic seed function and identify
the exceptional values $\lambda=\pm i$ for the inverse reconstruction
formula.}
\end{abstract}


\tableofcontents

\section{Introduction}

{\color{black}
The Weierstrass--Enneper representation is a cornerstone of minimal surface theory.
In the nineteenth century, it revealed a deep connection between complex analysis and differential geometry by providing explicit formulas for surfaces in $\mathbb{R}^3$ \cite{Nitsche1989, Osserman1986}.
Hoffman and Osserman extended this representation to $\mathbb{R}^4$ \cite{HO1980} and further developed it in a subsequent work \cite{HoffmanOsserman1983}.
 Jorge and Meeks studied global topological aspects of complete minimal surfaces of finite total curvature \cite{JM1983}.
More recent studies have presented new curvature identities and matrix-based formulations that deepen the understanding of this classical theory \cite{GYT2025}.
}

The Henneberg minimal surface was historically regarded as the only known example of a non-orientable minimal surface. However, subsequent developments in the literature have shown that this is not the case. In particular, Oliveira~\cite{DO1986} introduced new examples of non-orientable minimal surfaces, and later Martín~\cite{Ma2005} constructed complete non-orientable minimal surfaces in $\mathbb{R}^3$. Most recently, Toda and G\"{u}ler introduced a new version of the Weierstrass--Enneper representation in $\mathbb{R}^4$ \cite{TG2025}.

In this paper, we adapt the alternative Weierstrass--Enneper representation to construct explicit Henneberg-type minimal surfaces in $\mathbb{R}^4$. Our approach preserves the algebraic clarity of the classical formulas and remains equivalent to the holomorphic null curve method in $\mathbb{C}^4$. We present parametrizations, compute invariants, and provide projections that reveal the geometry of these surfaces.

\medskip 

{\color{black}
Our starting point is the holomorphic data $(f,g,h)$, which encodes a null curve in $\mathbb{C}^4$ whose real part integrates to a conformal minimal immersion in $\mathbb{R}^4$. We first recall how this scheme relates to the $(g_1,g_2)$ meromorphic and holomorphic one-form setting. We then transplant the classical Henneberg data into the $(f,g,h)$ framework and develop explicit parameterizations. Subsequent sections explain how the induced metric, the normals, and auxiliary vector fields arise directly from the formulas, followed by a discussion of the role of the parameter $\lambda$ and of higher-order families $\mathcal{H}_{(m,n)}$. 
}

\section{Henneberg-Type Minimal Surfaces in $\mathbb{R}^4$}

 {\color{black}
Weierstrass-type descriptions of minimal surfaces in $\mathbb{R}^4$ have a substantial earlier literature; in particular, Hoffman and Osserman's work on the generalized Gauss map is foundational~\cite{HO1980}
. We use the later work of Kawakami and Watanabe~\cite{KawakamiWatanabe2020} as a convenient reference for a two-meromorphic-function formulation in $\mathbb{R}^4$.
}

\subsection{The Alternative Representation versus the Kawakami--Watanabe Framework}

 {\color{black}
A standard meromorphic-function form of the Weierstrass representation for minimal surfaces in $\mathbb{R}^4$ is used, for example, in the work of Kawakami and Watanabe~\cite{KawakamiWatanabe2020}. In this formulation the immersion is encoded by two meromorphic functions $(g_1,g_2)$ together with a holomorphic one-form $h\,dz$, giving the holomorphic $1$-form 

}
\[
(\varphi_1,\varphi_2,\varphi_3,\varphi_4) 
= \Bigl( \dfrac{1}{2}(1+g_1g_2), 
\dfrac{\mathrm{i}}{2}(1-g_1g_2), 
\dfrac{1}{2}(g_1-g_2), 
-\dfrac{\mathrm{i}}{2}(g_1+g_2)\Bigr)\, h\,dz,
\]
which integrates to a minimal immersion into $\mathbb{R}^4$. 

{\color{black}
\begin{theorem}[Toda and G\"{u}ler~\cite{TG2025}]
Let $D\subset\mathbb C$ be a simply connected domain and let $f,g,h$ be holomorphic on $D$. Define
\begin{equation}\label{Weierstrass}
X(\omega, \bar{\omega}) = \operatorname{Re} \int
\bigg(
\dfrac{1}{2} f \left( 1 - g^2 - h^2 \right),\;
\dfrac{\mathrm{i}}{2} f \left( 1 + g^2 + h^2 \right),\;
fg,\;
fh
\bigg)\, d\omega.
\end{equation} 
 At every point at which $f(\omega)\neq0$, the map $X$ is a conformal minimal immersion. If $f$ has a zero, the Weierstrass differential vanishes there and the resulting map is a branched conformal minimal surface rather than a regular immersion at that point. On a non-simply-connected domain, the usual real-period conditions must additionally be imposed for $X$ to be globally single-valued.
\end{theorem}
}

{\color{black}
Our formulation above, by contrast, is expressed in terms of three holomorphic functions $(f,g,h)$. 
}

The Weierstrass differential in (\ref{Weierstrass}) is null, since $\sum_{k=1}^4\phi_k^2=0$.

 \textbf{Important Note.} The Kawakami--Watanabe representation is well-suited for studying Gauss map value distribution and 
Nevanlinna-type results, while the alternative $(f,g,h)$ representation has the advantage of a more direct analogy with the classical $\mathbb{R}^3$ case, and is more convenient for constructing explicit families of minimal surfaces 
by choosing simple polynomial data.

\subsection*{Notation and regularity conventions}
We use complex coordinates $\omega=u+\mathrm{i}v$ on a Riemann surface domain $D$. All derivatives are taken with respect to $u$ and $v$ unless stated otherwise. When writing $\Re\int\Phi$, the choice of base point affects only additive constants in $X$, irrelevant for intrinsic geometry. Regular points are those where $|f|\,\bigl(1+|g|^2+|h|^2\bigr)\neq 0$. Zeros of $f$ are branch points of the parametrization, whereas poles of the Weierstrass data correspond to punctures or ends of the parameter domain and require separate local analysis.

\section{Adaptation of the Henneberg Surface}
\subsection{Classical Henneberg Data in $\mathbb{R}^3$}

In the $(g,\omega)$ Weierstrass--Enneper convention for $\mathbb{R}^3$, 
a conformal minimal immersion is given by
\[
X(z)=\Re\int \left(\dfrac{1}{2}(1-g^2)\,\omega,\ 
\dfrac{\mathrm{i}}{2}(1+g^2)\,\omega,\ g\,\omega \right),
\]
with $g$ meromorphic and $\omega$ holomorphic. 
The classical \emph{Henneberg surface} corresponds to the data
\begin{equation}\label{eq:Henneberg3D}
g(z)=z, \qquad 
\omega=(1-z^{-4})\,dz, \qquad z\in\mathbb{C}\setminus\{0\}.
\end{equation}

{\color{black}
See Moya and P\'erez ~\cite[Eq.~(1) and Eq.~(5)]{Moya2023} and Osserman~\cite[Ch.~6]{Osserman1986}. 
}

Equivalently, in the $(g,dh)$ convention, one has $dh=g\,\omega=z(1-z^{-4})\,dz=(z-z^{-3})\,dz$. 

\medskip
\begin{remark}
The choice $\omega=(1-z^{-4})dz$ produces two types of terms after integration: a polynomial contribution encoding the interior geometry and a reciprocal part reflecting the effect of the puncture at $z=0$. Their interplay gives the well-known symmetries and saddle-type features of the Henneberg surface in projections.
\end{remark}

\subsection{Adaptation to $\mathbb{R}^4$}

To embed the Henneberg construction into our alternative Weierstrass--Enneper representation in $\mathbb{R}^4$ (see Toda and G\"{u}ler~\cite{TG2025}), we reinterpret the data 
(\ref{eq:Henneberg3D}) in terms of $(f,g,h)$. 

 {\color{black}
Throughout this paper, the term ``Henneberg-type'' is used in a precise constructional sense: the data retain the characteristic polynomial--reciprocal structure of the classical Henneberg Weierstrass data and reduce to known Henneberg constructions in distinguished subfamilies. We do not use the term to assert, for every $(m,n,\lambda)$, non-orientability, completeness, finite total curvature, or embeddedness. Those global properties require separate proofs. In particular, when $m=1$, $n$ is odd, and $\lambda=0$, the resulting three-dimensional data reproduce the corresponding odd member of the generalized Henneberg family studied by Moya and P\'erez~\cite{Moya2023}. Thus the novelty of the present construction lies in the $\mathbb R^4$ extension and in the additional complex parameter/data, rather than in rediscovering that three-dimensional subfamily.
}

 In $\mathbb{R}^4$, we adopt the following general Weierstrass data associated with the higher-order Henneberg type minimal surfaces family $\mathcal{H}_{(m,n)}$:
\begin{equation}\label{eq:HennebergR4data-general}
f(\omega)=2\,\omega^{-m-n-2}\,\bigl(\omega^{2m+2n}-1\bigr), \quad
g(\omega)=\omega^{m}, \quad
h(\omega)=\lambda\,\omega^{n}, 
\end{equation}
{\color{black}
where $\lambda\in\mathbb{C}$ and
$m,n\in\mathbb{N}_{\mathrm{odd}}$.
The restriction to positive odd integers preserves
the odd-order Henneberg-type subfamily and ensures
$m-n\neq\pm1$, thereby avoiding logarithmic terms
in the holomorphic primitives.
}

Integrating the associated holomorphic null differential yields the
following holomorphic primitive:

\[
\begin{aligned}
X_1(\omega) &= 
\frac{\omega^{\,m+n-1}}{m+n-1}
- \frac{\omega^{\,3m+n-1}}{3m+n-1}
- \lambda^{2}\frac{\omega^{\,m+3n-1}}{m+3n-1}
+ \frac{\omega^{-(m+n+1)}}{m+n+1} \\
&\quad + \frac{\omega^{\,m-n-1}}{m-n-1}
+ \lambda^{2}\frac{\omega^{\,n-m-1}}{n-m-1} + C_1, \\[1em]
X_2(\omega) &= \mathrm{i}\!\left(
\frac{\omega^{\,m+n-1}}{m+n-1}
+ \frac{\omega^{\,3m+n-1}}{3m+n-1}
+ \lambda^{2}\frac{\omega^{\,m+3n-1}}{m+3n-1}
+ \frac{\omega^{-(m+n+1)}}{m+n+1} \right. \\
&\quad \left.
- \frac{\omega^{\,m-n-1}}{m-n-1}
- \lambda^{2}\frac{\omega^{\,n-m-1}}{n-m-1}
\right)+C_2, \\[1em]
X_3(\omega) &= 2\left(
\frac{\omega^{\,2m+n-1}}{2m+n-1}
+ \frac{\omega^{-(n+1)}}{n+1}
\right) + C_3, \\[1em]
X_4(\omega) &= 2\lambda\left(
\frac{\omega^{\,m+2n-1}}{m+2n-1}
+ \frac{\omega^{-(m+1)}}{m+1}
\right)+C_4.
\end{aligned}
\]
{\color{black}
\medskip
The real part of the displayed holomorphic primitive defines the
Henneberg-type conformal minimal map on the regular set; zeros of $f$
are branch points. The expanded identity
$f(\omega)=2(\omega^{m+n-2}-\omega^{-(m+n+2)})$
is equivalent to~(\ref{eq:HennebergR4data-general}).

\noindent\textbf{Specialization $g(\omega)=\omega$, $h(\omega)=\lambda\omega$.}
For the same choice of 

}
\[
f(\omega)=2\big(\omega^{\,m+n-2}-\omega^{-(m+n+2)}\big),
\]
but with fixed
\[
g(\omega)=\omega,\qquad h(\omega)=\lambda\,\omega,\quad(\lambda\in\mathbb{C})
\]
the holomorphic $1$-form
\[
\Phi(\omega)
=\Bigl(\tfrac12 f(1-(1+\lambda^2)\omega^2),\ 
\tfrac{\mathrm i}{2}f(1+(1+\lambda^2)\omega^2),\ 
f\omega,\ 
f\lambda\omega\Bigr)\,d\omega
\]
integrates termwise to
\[
\int\!\Phi(\omega)
=\bigl(k_1(\omega),\,k_2(\omega),\,k_3(\omega),\,k_4(\omega)\bigr),
\]
where 
{\color{black}
\begin{align*}
k_1(\omega)
&=\frac{\omega^{m+n-1}}{m+n-1}
-\frac{(1+\lambda^2)\omega^{m+n+1}}{m+n+1}
+\frac{\omega^{-(m+n+1)}}{m+n+1}
-\frac{(1+\lambda^2)\omega^{-(m+n-1)}}{m+n-1},\\[4pt]
k_2(\omega)
&=\mathrm i\!\left(
\frac{\omega^{m+n-1}}{m+n-1}
+\frac{(1+\lambda^2)\omega^{m+n+1}}{m+n+1}
 + \frac{\omega^{-(m+n+1)}}{m+n+1}
 + \frac{(1+\lambda^2)\omega^{-(m+n-1)}}{m+n-1}\right),\\[4pt]
k_3(\omega)
&=\frac{2\omega^{m+n}}{m+n}
+\frac{2\omega^{-(m+n)}}{m+n},\\[4pt]
k_4(\omega)
&=\lambda\,k_3(\omega).
\end{align*}
}

\noindent
{\color{black}
\begin{remark}
The apparent sign difference between $X_2$ and $k_2$
arises from the different choices of Weierstrass data.
In particular, when $m=n=1$, the negative denominators
in $X_2$ reverse the corresponding signs, yielding
$X_2=k_2$ up to an additive constant.
Thus, the two expressions are consistent.
\end{remark}
}

{\color{black}
Finally, the immersion is given by
$
X(\omega)=\Re\bigl(k_1(\omega),\,k_2(\omega),\,k_3(\omega),\,k_4(\omega)\bigr).
$
\medskip
Since $f\,d\omega$ plays the role of $\omega$ in the $(g,h,\omega)$ formulation, we set on the punctured domain $D=\mathbb{C}\setminus\{0\}$ (by taking $m=n=1$) for Henneberg-type minimal surface $\mathcal{H}_{(1,1)}$:

}
\begin{equation}\label{eq:HennebergR4data-11}
f(\omega)=2(1-\omega^{-4}), \qquad 
g(\omega)=\omega, \qquad 
h(\omega)=\lambda\,\omega, \quad \lambda\in\mathbb{C}.
\end{equation}
 {\color{black}
For $\omega=u+\mathrm{i}v$, let
 
}
{\color{black}
\[
\phi (\omega)=
\left(
 \dfrac{1}{2} f (1-g^{2}-h^{2}),\;
 \dfrac{\mathrm{i}}{2} f (1+g^{2}+h^{2}),\;
 fg,\;
 fh
 \right)
 \]
denote the vector-valued holomorphic coefficient of the Weierstrass differential. Since $X=\Re\int\phi(\omega)\,d\omega$,
 
}
{\color{black}
\[
X_u=\Re\phi,\qquad X_v=-\Im\phi,\qquad X_u-\mathrm{i}X_v=\phi.
\]
}

{\color{black}
\medskip
\noindent\textit{Regularity and ends.}
For $m=n=1$, the data are holomorphic on $D=\mathbb C\setminus\{0\}$. The factor $1-\omega^{-4}$ in $f$ produces reciprocal terms near the puncture $\omega=0$ and polynomial terms at infinity. In particular,
}
{\color{black}
\[
fg=2(\omega-\omega^{-3}),\qquad
fh=2\lambda(\omega-\omega^{-3}),
\]
so the third and fourth holomorphic primitives have quadratic growth at infinity and reciprocal quadratic terms near the puncture. The first two components contain cubic polynomial terms and reciprocal cubic terms. The zeros of $f$, namely $\omega^4=1$, are branch points of the parametrization. Accordingly, the map is a regular conformal minimal immersion on
$\mathbb C\setminus(\{0\}\cup\{\omega:\omega^4=1\})$ and extends locally across the four zeros as a branched conformal minimal surface.
}

\subsection{Case study: values of the parameter $\lambda$}

\subsubsection*{Case 1: $\lambda\in \mathbb{C}$.}
When $m=n=1$, $\lambda \in \mathbb{C}$, Henneberg type minimal curve in complex 4-space is given by 
{\color{black}
\[
X(\omega)=
 \begin{pmatrix}
\omega - \dfrac{1+\lambda^{2}}{3}\,\omega^{3} + \dfrac{1}{3}\,\omega^{-3} - (1+\lambda^{2})\,\omega^{-1} \\[0.9em]
\mathrm{i}\!\left( \omega + \dfrac{1+\lambda^{2}}{3}\,\omega^{3} + \dfrac{1}{3}\,\omega^{-3} + (1+\lambda^{2})\,\omega^{-1} \right) \\[0.9em]
\omega^{2}+\omega^{-2} \\[0.9em]
\lambda\!\left(\omega^{2}+\omega^{-2}\right)
\end{pmatrix}.
\]
}

In order to obtain an explicit parametrization in real coordinates, 
we expand all components of the immersion $X(u,v)$ 
in terms of the parameters $u,v \in \mathbb{R}$ and the complex constant 
$\lambda = \alpha + \mathrm{i}\beta$, where $\alpha = \Re(\lambda)$ and $\beta = \Im(\lambda)$. The resulting immersion $\mathcal{H}_{(1,1)}$ in $\mathbb{R}^4$ takes the form 

{\color{black}
\[
 \begin{aligned}
x(u,v) &= u
 - \frac{1+\alpha^{2}-\beta^{2}}{3}\,(u^{3}-3uv^{2})
 - \frac{1+\alpha^{2}-\beta^{2}}{u^{2}+v^{2}}\,u
+ \frac{u^{3}-3uv^{2}}{3(u^{2}+v^{2})^{3}}
\\
&\quad
 + \frac{2\alpha\beta}{3}\,(3u^{2}v - v^{3})
 \textcolor{black}{- \frac{2\alpha\beta}{u^{2}+v^{2}}\,v},
\\[1.3em]
y(u,v) &= -v
 - \frac{1+\alpha^{2}-\beta^{2}}{3}\,(3u^{2}v - v^{3})
 + \frac{1+\alpha^{2}-\beta^{2}}{u^{2}+v^{2}}\,v
 + \frac{3u^{2}v - v^{3}}{3(u^{2}+v^{2})^{3}}\\
&\quad
 - \frac{2\alpha\beta}{3}\,(u^{3}-3uv^{2})
 - \frac{2\alpha\beta}{u^{2}+v^{2}}\,u,
\\[1.3em]
z(u,v) &= (u^{2}-v^{2})\!\left(1 + \frac{1}{(u^{2}+v^{2})^{2}}\right),
\\[1.3em]
w(u,v) &=
\alpha (u^{2}-v^{2})\!\left(1+\frac{1}{(u^{2}+v^{2})^{2}}\right)
 - 2\beta\, uv\!\left(1-\frac{1}{(u^{2}+v^{2})^{2}}\right),
\end{aligned} 
\]
}
where
\[
X(u,v)=(x(u,v),\, y(u,v),\, z(u,v),\, w(u,v)).
\]

{\color{black}
\noindent\textit{Asymptotics and projections.}
The large $|\omega|$ behavior is dominated by the cubic and quadratic monomials visible in the first two and last two coordinates, respectively, whereas the reciprocal terms control the geometry near $\omega=0$. Orthogonal projections onto three-dimensional coordinate subspaces are used here only as visualizations. Such projections do not preserve conformality in general: if $Y=(x_1,x_2,x_3)$ is obtained from a conformal immersion $X=(x_1,x_2,x_3,x_4)$, then $|Y_u|^2=|X_u|^2-x_{4,u}^2$, $|Y_v|^2=|X_v|^2-x_{4,v}^2$, and $\langle Y_u,Y_v\rangle=-x_{4,u}x_{4,v}$. Hence no conformality claim is made for three-dimensional coordinate projections.
}

 \subsubsection*{Case 2: $\lambda\in\mathbb{R}$.}

 {\color{black}
For $m=n=1$ and real $\lambda$, the tangent geometry is obtained directly from the holomorphic null differential. With

}
{\color{black}
\[
f=2(1-\omega^{-4}),\qquad g=\omega,\qquad h=\lambda\omega,
\]
}
 {\color{black}
set
}
{\color{black}
\[
\phi=
\left(
\frac{1}{2} f(1-(1+\lambda^2)\omega^2),\;
\frac{\mathrm i}{2}f(1+(1+\lambda^2)\omega^2),\;
f\omega,\;
\lambda f\omega
\right).
\]
}
{\color{black}
Then
}
{\color{black}
\[
X_u=\Re\phi,\qquad X_v=-\Im\phi.
\]
}

{\color{black}
The nullity relation $\sum_{j=1}^4\phi_j^2=0$ yields

}
{\color{black}
\[
\langle X_u,X_v\rangle=0,\qquad
|X_u|^2=|X_v|^2=E,
\]
}
where
{\color{black}
\[
E=
\frac{|f|^2}{4}
\left(1+(1+\lambda^2)|\omega|^2\right)^2.
\]
}

{\color{black}
For arbitrary complex $\lambda$, the corresponding identity is

}
{\color{black}
\[
E=\frac{|f|^2}{4}
\left(
1+2(1+|\lambda|^2)|\omega|^2
+|1+\lambda^2|^2|\omega|^4
\right).
\]
}

{\color{black}
Hence, at every regular point,

}
{\color{black}
\[
\mathbf e_1=\frac{X_u}{\sqrt E},\qquad
\mathbf e_2=\frac{X_v}{\sqrt E}
\]
is an orthonormal tangent frame. A local orthonormal normal frame can be constructed by normal projection followed by Gram--Schmidt. For a constant vector $a\in\mathbb R^4$, define

}
{\color{black}
\[
P_N(a)=
a-\frac{\langle a,X_u\rangle}{E}X_u
-\frac{\langle a,X_v\rangle}{E}X_v.
\]
Choose constant vectors $a,b\in\mathbb R^4$ locally so that the following denominators do not vanish. Then
}
{\color{black}
\[
\mathbf n_1=\frac{P_N(a)}{|P_N(a)|},
 \]
}
{\color{black}
\[
\widetilde{\mathbf n}_2
=P_N(b)-\langle P_N(b),\mathbf n_1\rangle\mathbf n_1,
\qquad
\mathbf n_2=
\frac{\widetilde{\mathbf n}_2}{|\widetilde{\mathbf n}_2|}.
 \]
Thus
$\{\mathbf e_1,\mathbf e_2,\mathbf n_1,\mathbf n_2\}$
is a local orthonormal frame along the regular part of the surface.
}
{\color{black}
For real $\lambda$, the corresponding parametrization of $\mathcal H_{(1,1)}$ is
}
{\color{black}
\[
X(u,v) =
 \begin{pmatrix}
u - \dfrac{1+\lambda^2}{3}\,(u^3-3uv^2) 
 + \dfrac{u^3-3uv^2}{3(u^2+v^2)^3} 
 - (1+\lambda^2)\dfrac{u}{u^2+v^2} \\[1.2em]
-\,v - \dfrac{1+\lambda^2}{3}\,(3u^2v-v^3) 
 + \dfrac{3u^2v-v^3}{3(u^2+v^2)^3} 
 + (1+\lambda^2)\dfrac{v}{u^2+v^2} \\[1.2em]
u^2-v^2 + \dfrac{u^2-v^2}{(u^2+v^2)^2} \\[1.2em]
\lambda\!\left(u^2-v^2+\dfrac{u^2-v^2}{(u^2+v^2)^2}\right)
\end{pmatrix}.
\]
}

\subsubsection*{Case 2.1: $\lambda=0$.}
When $\lambda=0$, the data reduce to
\[
f(\omega)=2(1-\omega^{-4}), \quad g(\omega)=\omega, \quad h(\omega)=0.
\]
{\color{black}
The immersion \(X(u,v)\) then has a trivial fourth coordinate, and the first three coordinates coincide with the classical Henneberg surface in $\mathbb{R}^3$. 
Thus $\lambda=0$ recovers the classical Henneberg construction in the coordinate subspace $\mathbb{R}^3\times\{0\}\subset\mathbb{R}^4$.
}

\subsubsection*{Case 2.2: $\lambda\in\mathbb{R}\setminus\{0\}$.}
For $m=n=1$ and real $\lambda\neq0$, the identity $fh=\lambda fg$ implies $w=\lambda z$ (up to the same choice of additive constants). Consequently the image is contained in the three-dimensional linear subspace
{\color{black}
\[
\{(x,y,z,w)\in\mathbb R^4: w=\lambda z\}.
\]
Thus this real-$\lambda$ subfamily does \emph{not} acquire a genuinely four-dimensional image. The parameter $\lambda$ changes the first two coordinate functions through the factor $1+\lambda^2$ and changes the position of the containing three-space, while the image remains three-dimensional for $m=n=1$.
}

{\color{black}
\begin{example}
Let $m=n=1$ and $\lambda=1+i$. Therefore, the real parametric form of the minimal surface $\mathcal{H}_{(1,1)}\subset\mathbb{R}^4$ is determined by 
 \[
 \begin{aligned}
x(u,v) &=
 u - \dfrac{1}{3}(u^3 - 3uv^2) + \dfrac{2}{3}(3u^2v - v^3)
 - \dfrac{u+2v}{u^2+v^2}
 + \dfrac{u^3 - 3uv^2}{3(u^2+v^2)^3}
 \\[12pt]
 y(u,v) &= 
 - v - \dfrac{1}{3}(3u^2v - v^3) - \dfrac{2}{3}(u^3 - 3uv^2)
+ \dfrac{v-2u}{u^2+v^2}
 + \dfrac{3u^2v - v^3}{3(u^2+v^2)^3}
 \\[12pt]
 z(u,v) &= 
 (u^2 - v^2)\left(1 + \dfrac{1}{(u^2+v^2)^2}\right)
 \\[12pt]
 w(u,v) &= 
 (u^2 - v^2)\left(1 + \dfrac{1}{(u^2+v^2)^2}\right)
- 2uv\left(1 - \dfrac{1}{(u^2+v^2)^2}\right).
 \end{aligned} 
 \]
The polar form of the minimal surface $\mathcal{H}_{(1,1)}$ is given by 
 \[
\begin{aligned}
x(r,\theta)
&= \left(r - \frac{1}{r}\right)\cos\theta
- \frac{2}{r}\sin\theta
- \frac{1}{3}\left(r^3 - \frac{1}{r^3}\right)\cos(3\theta)
+ \frac{2}{3}r^3\sin(3\theta),
\\[10pt]
y(r,\theta) &=
-\left(r - \frac{1}{r}\right)\sin\theta
- \frac{2}{r}\cos\theta
- \frac{1}{3}\left(r^3 - \frac{1}{r^3}\right)\sin(3\theta)
- \frac{2}{3}r^3\cos(3\theta),
\\[10pt]
z(r,\theta) &=
\left(r^2 + \frac{1}{r^2}\right)\cos(2\theta),
\\[10pt]
w(r,\theta) &=
\left(r^2 + \frac{1}{r^2}\right)\cos(2\theta)
- \left(r^2 - \frac{1}{r^2}\right)\sin(2\theta).
\end{aligned}
\] 
\end{example}
}

\begin{figure}[t!]
 \centering
 \IfFileExists{1.jpg}{\includegraphics[width=1\linewidth]{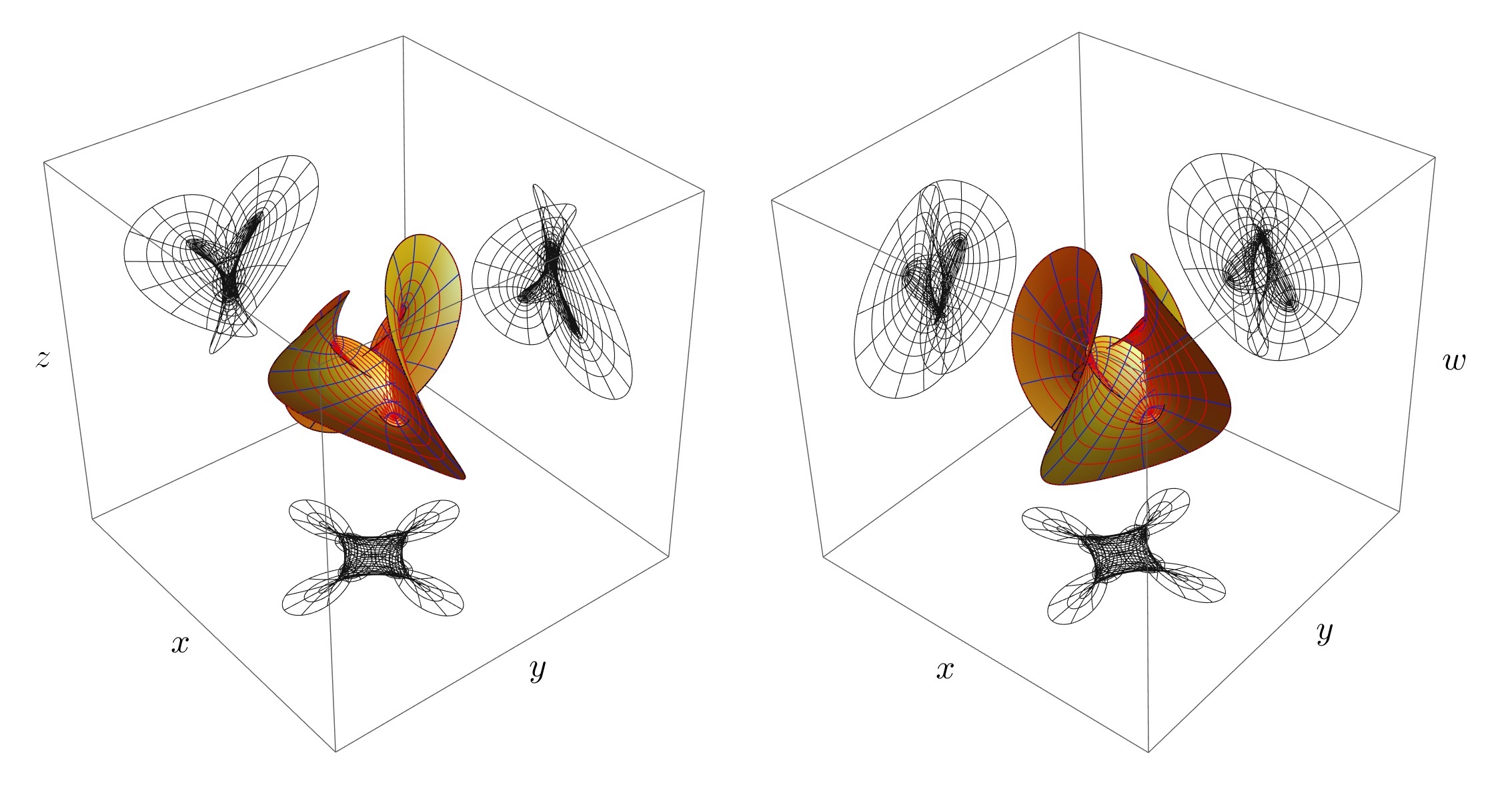}}{\fbox{\parbox[c][35mm][c]{0.85\linewidth}{\centering Figure 1: insert the Mathematica image \texttt{1.jpg}.}}}
 \caption{Coordinate projections of the $\mathbb{R}^4$ minimal surface $\mathcal{H}_{(1,1)}$ in $xyz$-space (left) and $xyw$-space (right), where $\lambda = 1 + \mathrm{i}$.}
 \label{fig:xyz-xyw-proj-a}
\end{figure}

\begin{figure}[t!]
 \centering
 \IfFileExists{2.jpg}{\includegraphics[width=1\linewidth]{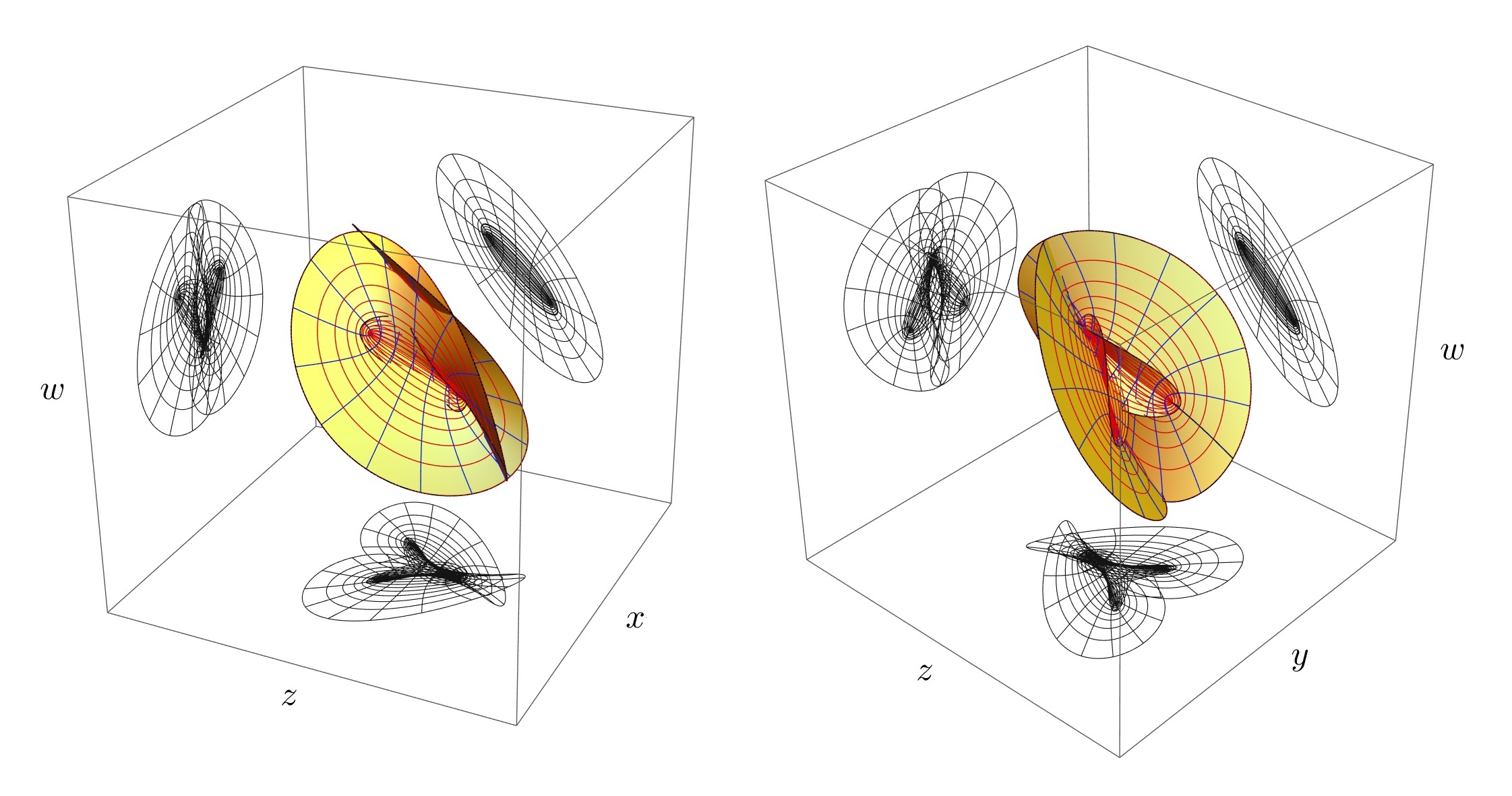}}{\fbox{\parbox[c][35mm][c]{0.85\linewidth}{\centering Figure 2: insert the Mathematica image \texttt{2.jpg}.}}}
 \caption{Coordinate projections of the $\mathbb{R}^4$ minimal surface $\mathcal{H}_{(1,1)}$ in $xzw$-space (left) and $yzw$-space (right), where $\lambda = 1 + \mathrm{i}$.}
 \label{fig:xzw-yzw-proj-a}
\end{figure}

Figures~\ref{fig:xyz-xyw-proj-a} and~\ref{fig:xzw-yzw-proj-a} illustrate the
coordinate projection surfaces obtained from the
polar parametrization of the minimal surface
$\mathcal{H}_{(1,1)}$ in Example~1.
These projected surfaces are not necessarily minimal.

{\color{black}
\begin{proposition}[Gaussian curvature of $\mathcal{H}_{(1,1)}$]
Let
\[
a=1+|\lambda|^2,\qquad
b=|1+\lambda^2|^2,\qquad
P(r)=1+2ar^2+br^4,
\]
where $r=|\omega|$. At every regular point,
the induced metric is
\[
ds^2=\frac{|f(\omega)|^2}{4}P(r)(du^2+dv^2),
\qquad f(\omega)=2(1-\omega^{-4}),
\]
and the Gaussian curvature is
\[
K=-\frac{16(a+2br^2+abr^4)}
{|f(\omega)|^2P(r)^3}.
\]
In particular, $K<0$ at every regular point.
\end{proposition}

\begin{proof}
Since
\[
E=G=\frac{|f|^2}{4}P(r),\qquad F=0,
\]
and $\log|f|^2$ is harmonic away from
the zeros and poles of $f$, we obtain
\[
K=-\frac{1}{2E}\Delta\log P,
\qquad
\Delta\log P=
\frac{8(a+2br^2+abr^4)}{P(r)^2}.
\]
The result follows by substitution.
\end{proof}
}

\section{Integral-Free Representation in $\mathbb{R}^4$}

 {\color{black}
In this section, $\mathcal F$ denotes the scalar holomorphic seed function, whereas $\Phi$ denotes the vector-valued Weierstrass differential introduced above.
}

 {\color{black}
The Weierstrass--Enneper formalism represents minimal surfaces using holomorphic data 
$(f,g)$ in $\mathbb{R}^3$ and $(f,g,h)$ in $\mathbb{R}^4$~\cite{Nitsche1989,TG2025}. 
The \emph{integral-free} form derives the immersion from a single holomorphic 
seed function $\mathcal{F}$:
 
}
{\color{black}
\[
f = \mathcal{F} ''', \qquad g(\omega) = \omega, \qquad h(\omega) = \lambda\,\omega, 
\quad \lambda \in \mathbb{C}.
\]
}

{\color{black}
Let $\mathcal{F}$ be holomorphic on a simply connected domain $D\subset\mathbb{C}$. Define

}
\[
X(\omega)=\Re\,\bigl(k_1(\omega),\,k_2(\omega),\,k_3(\omega),\,k_4(\omega)\bigr),
\]
where
{\color{black}
\[
 \begin{aligned}
k_1 &= \tfrac12(1-(1+\lambda^2)\omega^2)\mathcal{F}''+(1+\lambda^2)\omega\mathcal{F}'-(1+\lambda^2)\mathcal{F},\\
k_2 &= \tfrac{\mathrm i}{2}(1+(1+\lambda^2)\omega^2)\mathcal{F}''-\mathrm i(1+\lambda^2)\omega\mathcal{F}'+\mathrm i(1+\lambda^2)\mathcal{F},\\
k_3 &= \omega\mathcal{F}''-\mathcal{F}',\qquad
k_4 = \lambda(\omega\mathcal{F}''-\mathcal{F}').
\end{aligned} 
\]
}

 {\color{black}
For $1+\lambda^2\neq0$ (equivalently, $\lambda\neq\pm i$), the inverse reconstruction identity is
 
}
{\color{black}
\[
 \mathcal{F} =\frac{((1+\lambda^{2})\omega^{2}-1)}{2(1+\lambda^{2})}k_{1}
- \mathrm i\,\frac{((1+\lambda^{2})\omega^{2}+1)}{2(1+\lambda^{2})}k_{2}
- \frac{\omega}{(1+\lambda^{2})}(k_{3}+\lambda k_{4}).
\]
}

{\color{black}
For $\lambda=\pm i$, this inverse formula is singular and is not used; the forward formulas for $k_1,\ldots,k_4$ remain meaningful. For $\lambda=0$, $X(\omega)$ reduces to the integral-free Weierstrass form in $\mathbb{R}^3$.
}

{\color{black}
For the specialized data $g(\omega)=\omega$ and $h(\omega)=\lambda\omega$, the seed function corresponding to $f(\omega)=2(\omega^{m+n-2}-\omega^{-(m+n+2)})$ is
}
{\color{black}
\[
 \mathcal{F} (\omega)
= \frac{2}{(m+n-1)(m+n)(m+n+1)}
\left(\omega^{m+n+1} + \omega^{-m-n+1}\right).
\]
}
In particular, for $\mathcal{H}_{(1,1)}$ the function simplifies to
{\color{black}
\[
 \mathcal{F} (\omega) = \frac{1}{3}\,(\omega^{3} + \omega^{-1}).
\]
}

\section{Discussion and Further Questions}
The construction developed here is an explicit Weierstrass-data construction in $\mathbb R^4$. The designation `` Henneberg-type '' refers to the inherited polynomial--reciprocal structure and to the reduction to classical/generalized Henneberg data in the subfamilies identified above; it is not intended as a blanket assertion of the global properties of the classical Henneberg surface. The regular/branched distinction is explicit: zeros of $f$ are branch points, and all local immersion statements are understood away from those zeros.

{\color{black}
Natural further questions concern which subfamilies are complete, have finite total curvature, or admit a non-orientable quotient via an antiholomorphic involution of the oriented double cover. These are global questions and are separated from the local explicit construction proved in the present paper. A systematic treatment of them will require period, end, and topology analyses beyond the algebraic parametrizations established here.
}

\vspace{2em}
\noindent

\textbf{Erhan Güler}\\
Department of Mathematics and Statistics, Texas Tech University,\\
Lubbock, TX 79409, USA\\
\textit{Email:} \texttt{eguler@ttu.edu}\\[1em]

\textbf{Magdalena Toda}\\
Department of Mathematics and Statistics, Texas Tech University,\\
Lubbock, TX 79409, USA\\
\textit{Email:} \texttt{magda.toda@ttu.edu}


\begin{thebibliography}{99}



\bibitem{DO1986} 
 M.E.G.G. De Oliveira, Some new examples of nonorientable minimal surfaces, \textit{Proc. Amer. Math. Soc.}, \textbf{98}(4) (1986), 629--636.

\bibitem{GYT2025}
E. G\"{u}ler, Y. Yayl{\i}, M. Toda, Differential geometry and matrix-based generalizations of the Pythagorean theorem in space forms, \textit{Mathematics}, \textbf{13}(5) (2025), 836.

\bibitem{HO1980} 
D.A. Hoffman and R. Osserman, \textit{The geometry of the generalized Gauss map}, Mem. Amer. Math. Soc. 28, no. 236, 1980.

\bibitem{HoffmanOsserman1983}
D.A. Hoffman and R.~Osserman,
\newblock The Gauss map of surfaces in $\mathbb{R}^n$,
\newblock \emph{J. Differential Geom.}, 18: (1983) 733--754.

\bibitem{JM1983} 
L.P. Jorge and W.H. Meeks, The topology of complete minimal surfaces of finite Gaussian curvature, \textit{Topology}, \textbf{22} (1983), 203--221.

\bibitem{KawakamiWatanabe2020}
Y.~Kawakami and M.~Watanabe, 
\newblock The Gauss images of complete minimal surfaces of genus zero of finite total curvature,
\emph{J. Geom. Anal.}, 34: 270 (2024).

\bibitem{Ma2005}
F. Martín,
Complete nonorientable minimal surfaces in $\mathbb{R}^3$,
in \emph{Global theory of minimal surfaces},
Clay Math. Proc., 2, Amer. Math. Soc., Providence, RI, 371--380 (2005).

\bibitem{Moya2023}
D.~Moya and J. Perez, Generalized Henneberg stable minimal surfaces, \emph{Results Math.} 78:53 (2023), 1--25.

\bibitem{Nitsche1989}
J.C.C. Nitsche, \textit{Lectures on minimal surfaces}, Cambridge University Press, 1989.

\bibitem{Osserman1986}
R. Osserman, \textit{A survey of minimal surfaces}, Dover Publications, 1986.

\bibitem{TG2025}
M. Toda and E. G\"{u}ler, Generalized Weierstrass--Enneper representation for minimal surfaces in $\mathbb{R}^4$, \emph{AIMS Mathematics}, 10(9): 22406--22420 (2025). 


\end{thebibliography}
\end{document}